\documentclass[12pt]{article}
\usepackage{graphicx}
\usepackage{amssymb}
\usepackage{amsmath}
\numberwithin{equation}{section}

\begin{document}

\author{Lev Sakhnovich}
\date{November,24, 2005}
\textbf{On Sample Functions Behavior of Stable Processes}
\begin{center}\textbf{Lev Sakhnovich}\end{center}
\emph{735 Crawford ave., Brooklyn, 11223, New York, USA}\\
 E-mail:Lev.Sakhnovich@verizon.net\\
\begin{center}\textbf{Abstract}\end{center}
We investigate the asymptotic behavior  of sample functions of
stable processes when $t{\to}\infty$. We compare our results with
the iterated logarithm law, results for the first hitting time and
most visited sites problems.\\
\section{Main Notions}
1. Let $X_{1},X_{2}...$ be mutually independent random variables
with the same law of distribution $F(x)$. The distribution $F(x)$
is called \emph{strictly stable} if the random variable
\begin{equation}
X=(X_{1}+X_{2}+...+X_{n})/n^{1/\alpha} \end{equation} is also
distributed according to the law $F(x)$. The number $\alpha \quad
(0<\alpha{\leq}2)$ is called a \emph{characteristic exponent} of
the distribution. When $\alpha{\ne}1$ the characteristic function
of $X$ has the form (see[2])
\begin{equation}E[\exp{(i{\xi}X)}]=\exp{\{-\lambda|\xi|^{\alpha}[1-i\beta({\mathrm{sign}}{\xi})(\tan{\frac{\pi\alpha}{2}})]\}},
\end{equation} where $-1{\leq}\beta{\leq}1,\quad {\lambda}>0$.
When $\alpha=1$ we have
\begin{equation}E[\exp{(i{\xi}X)}]=\exp{\{-\lambda|\xi|[1+\frac{2i\beta}{\pi}({\mathrm{sign}}{\xi})(\log{|\xi|})]\}},
\end{equation} where $-1{\leq}\beta{\leq}1,\quad {\lambda}>0$. The
homogeneous process $X(\tau) (X(0)=0)$ with independent increments
is called a stable process if
\begin{equation}E[\exp{(i{\xi}X(\tau))}]=\exp{\{-\tau|\xi|^{\alpha}[1-i\beta({{\mathrm{sign}}\xi})(\tan{\frac{\pi\alpha}{2}})]\}},
\end{equation} where $0<\alpha<2, \alpha{\ne}1,-1{\leq}\beta{\leq}1,\quad {\tau}>0$.
When $\alpha=1$ we have
\begin{equation}E[\exp{(i{\xi}X(\tau))}]=\exp{\{-\tau|\xi|[1+\frac{2i\beta}{\pi}({{\mathrm{sign}}\xi})(\log{|\xi|)}]\}},
\end{equation} where $-1{\leq}\beta{\leq}1,\quad {\tau}>0$. The
stable processes are a natural generalization of the Wiener
processes. In many theoretical and applied problems it is
important to estimate the value \begin{equation}
p_{\alpha}(t,a)=P(\sup{|X(\tau)|}<a),\quad
0{\leq}\tau{\leq}t.\end{equation}It is proved (see [6],[7],[9])
that the value of $p_{\alpha}(t,a)$ decreases very quickly by the
exponential law when $t{\to}\infty$. This fact prompted the idea
to consider the case when the value of $a$ in (1.6) depends on $t$
and $a(t){\to}\infty,\quad t{\to}\infty$. In this paper we deduce
the conditions under which one of the following three cases is
realized:\\ 1) $\lim{p_{\alpha}(t,a(t))}=1,\quad, t{\to}\infty.$
\\ 2) $\lim{p_{\alpha}(t,a(t))}=0, \quad t{\to}\infty.$\\ 3)
$\lim{ p_{\alpha}(t,a(t))}=p_{\infty}, \quad  0<p_{\infty}{\leq}1,
\quad t{\to}\infty.$\\ The value of $p_{\alpha}(t,a(t))$ is the
probability of the trajectory of the process  $X(\tau)$ that
remains
 inside the corridor $|X(\tau)|<a(t)$
 when
$0{\leq}\tau{\leq}t$.\\ We investigate the situation when
$t{\to}0$ too. We consider separately the classical case when
$\alpha=2$. We think that even in this case our results are new.
In the last part of the paper we compare our results with the
classical iterated logarithm law (see [10]) , the first hitting
time problem (see[2],[5],[10]) and the most visited sites results
(see [1]).\\ \textbf{Remark 1.1.} In the famous work by M.Kac [3]
the connection of the theory of stable processes and the theory of
integral equations was shown. M.Kac considered in detail only the
case $\alpha=1,\quad \beta=0.$  The case $0<\alpha<2,\quad
\beta=0$ was later studied by H.Widom [11]. As to the general case
$0<\alpha<2,\quad -1{\leq}\beta{\leq}1 $ it was investigated in
our works [6],[7],[9]. In all the mentioned works the parameter
$a$ in (1.6) was fixed. The present paper is dedicated to the
important case when $a$ depends on $t$ and $a(t){\to}\infty,\quad
t{\to}\infty$.
\section{Auxillary results}
1. In this section we  formulate some results from our paper [7]
(see also [9], Ch.7). An important role in our approach is played
by formula [3] \begin{equation}
\int_{0}^{\infty}e^{-su}p_{\alpha}(u,a)du=\int_{-a}^{a}\psi_{\alpha}(x,s,a)dx.\end{equation}
Here $\psi_{\alpha}(x,s,a)$ is defined by relation
\begin{equation}
\psi_{\alpha}(x,s,a)=(I+sB_{\alpha}^{\star})^{-1}\Phi_{\alpha}(0,x,a),\end{equation}
(The operator $B_{\alpha}$ and its kernel $\Phi_{\alpha}(x,y,a)$
will be introduced later.\\Further we consider the three cases.\\
\emph{Case 1.} $0<\alpha<2,\quad \alpha{\ne}1,\quad -1<\beta
<1.$\\ \emph{Case 2.} $1<\alpha<2,\quad  \beta={\pm}1.$\\
\emph{Case 3.} $\alpha=1,\quad \beta=0.$\\ Now we introduce the
operators \begin{equation} B_{\alpha}f=
\int_{-a}^{a}\Phi_{\alpha}(x,y,a)f(y)dy \end{equation} acting in
the space $L^{2}(-a,a)$. \\In \emph{case 1} the kernel
$\Phi_{\alpha}(x,y,a)$ has the following form
(see[7],[9])\begin{equation}
\Phi_{\alpha}(x,y,a)=C_{\alpha}(2a)^{\mu-1}\int_{a|x-y|}^{a^{2}-xy}{[z^{2}-a^{2}(x-y)^{2}]}^{-\rho}[z-a(x-y)]^{2\rho-\mu}dz,\end{equation}
where the constants $\mu, \rho,$ and $C_{\alpha}$ are defined by
the relations $\mu=2-\alpha,$ \begin{equation}
\sin{\pi\rho}=\frac{1-\beta}{1+\beta}\sin{\pi(\mu-\rho)},\quad
0<\mu-\rho<1,\end{equation}\begin{equation}C_{\alpha}=
\frac{\sin\pi{\rho}}{(\sin{\pi\alpha/2})(1-\beta)\Gamma(1-\rho)\Gamma(1+\rho-\mu)}.\end{equation}
Here $\Gamma(z)$ is Euler's gamma function. We remark that the
constants $\mu, \rho,$ and $C_{\alpha}$ do not depend on parameter
$a$.\\ In \emph{case 2} when $\beta=1$ we have [7],[9]
\begin{equation} \Phi_{\alpha}(x,y,a)=\frac
{(\cos{\pi\alpha/2})}{(2a)^{\alpha
-1}\Gamma(\alpha)}\{[a|x-y|+y-x]^{\alpha-1}-(a-x)^{\alpha-1}(a+y)^{\alpha-1}\}
\end{equation}
In \emph{case 2} when $\beta=-1$ we have [7],[9]
\begin{equation} \Phi_{\alpha}(x,y,a)=\frac
{(\cos{\pi\alpha/2})}{(2a)^{\alpha-1}\Gamma(\alpha)}\{[a|x-y|+x-y]^{\alpha-1}-(a+x)^{\alpha-1}(a-y)^{\alpha-1}\}
\end{equation}
Finally, in \emph{case 3} according to M.Kac [3] we have
\begin{equation}
\Phi_{1}(x,y,a)=\frac{1}{4}{\mathrm{log}}\frac{[a^{2}-xy+\sqrt{(a^{2}-x^{2})(a^{2}-y^{2})}]}{[a^{2}-xy-\sqrt{(a^{2}-x^{2})(a^{2}-y^{2})}]}\end{equation}
\section{The asymptotic behavior of $p_{\alpha}(t,a)$}
1. In this section we use the  assertion, which immediately
follows from relations (2.4),(2.7)-(2.9).\\ \textbf{Proposition
3.1.} \emph{In all the three cases the relation}
\begin{equation}
\Phi_{\alpha}(ax,ay,a)=a^{\alpha-1}\Phi_{\alpha}(x,y,1)
\end{equation} \emph{is valid}.\\
From equalities (2.1), (2.2), and (3.1) we deduce the assertion.\\
\textbf{Corollary 3.1.} \emph{In all the three cases the following
relations}
\begin{equation}
\psi_{\alpha}(ax,s,a)=a^{\alpha-1}\psi_{\alpha}(x,sa^{\alpha},1),\end{equation}
\begin{equation}p_{\alpha}(t,a)=p_{\alpha}(\frac{t}{a^{\alpha}},1)\end{equation}
\emph{hold}.\\ 2. The kernel $\Phi_{\alpha}(x,y,a)$ of the
operator $B_{\alpha}$ is non-negative [7],[9]. So the operator
$B_{\alpha}$ has the  non-negative eigenfunction
$g_{\alpha}(x,a)$. The corresponding eigenvalue
$\lambda_{\alpha}(a)$ is positive and larger than the modules of
any other eigenvalues of the operator $B_{\alpha}$. The eigenvalue
of the adjoint operator $B_{\alpha}^{\star}$ with  the largest
modulus is also $\lambda_{\alpha}(a)$. We denote the corresponding
non-negative eigenfunction by $h_{\alpha}(x,a)$. We normalize the
functions $g_{\alpha}(x,a)$ and $h_{\alpha}(x,a)$ by the condition
$(g_{\alpha},h_{\alpha})=1$ . Now we can formulate the following
result.\\ \textbf{Theorem 3.1.} \emph{Let one of the following
conditions be fulfilled:}\\ I. $0<\alpha<2,\quad
\alpha{\ne}1,\quad  -1<\beta<1.$\\ II. $1<\alpha<2,\quad
\beta={\pm}1.$\\ III. $\alpha=1, \quad \beta=0.$\\ \emph{Then the
asymptotic equality holds} \begin{equation}
p_{\alpha}(t,a)=e^{-t/\lambda_{\alpha}(a)}g_{\alpha}(0,a)\int_{-a}^{a}h_{\alpha}(x,a)dx[1+o(1)],
\quad t{\to}\infty.\end{equation}We introduce the notations
\begin {equation}\lambda_{\alpha}(1)=\lambda_{\alpha},\quad
p_{\alpha}(t,1)=p_{\alpha}(t),\quad
g_{\alpha}(x,1)=g_{\alpha}(x),\quad
h_{\alpha}(x,1)=h_{\alpha}(x).\end{equation}Using relation (3.3)
and notations (3.5) we can rewrite formula (3.4) in the following
way
\begin{equation}
p_{\alpha}(t,a)=e^{-t/[a^{\alpha}\lambda_{\alpha}]}g_{\alpha}(0)\int_{-1}^{1}h_{\alpha}(x)dx[1+o(1)],
\quad t{\to}\infty.\end{equation}\textbf{Remark 3.1.} The operator
$B_{\alpha}$ is self-adjoint when $\beta=0$. In this case
$h_{\alpha}=g_{\alpha}$.\\ \textbf{Remark 3.2.} The value
$\lambda_{\alpha}$ characterizes how fast $p_{\alpha}(t,a)$
converges to zero when $t{\to}\infty$. The two-sided estimation
for $\lambda_{\alpha}$ when $\beta=0$ is given in [4] (see also
[9],p.150).\\ 3. Now we consider the case when the parameter $a$
depends on $t$. From Theorem 3.1.  we deduce the assertion.\\
\textbf{Corollary 3.2.} \emph{Let one of conditions I-III of
Theorem 3.1 be fulfilled and}\begin{equation}
\frac{t}{a^{\alpha}(t)}{\to}\infty,\quad
t{\to}\infty.\end{equation}\emph{Then the following equalities are
true:}\begin{equation} 1)\quad
p_{\alpha}(t,a(t))=e^{-t/[a^{\alpha}(t)\lambda_{\alpha}]}g_{\alpha}(0)\int_{-1}^{1}h_{\alpha}(x)dx[1+o(1)],
\quad t{\to}\infty.\end{equation}
\begin{equation}
2)\quad \lim{p_{\alpha}(t,a)}=0,\quad t{\to}\infty.\end{equation}
\begin{equation}3)\quad \lim{ P[\sup{|X(\tau)|}}>a(t)]=1,\quad
0{\leq}\tau{\leq}t,\quad t{\to}\infty. \end{equation}
\textbf{Remark 3.3.} Condition (3.7) is equivalent to the
condition
\begin{equation} \frac {a(t)}{t^{1/\alpha}}{\to}0,\quad t{\to}\infty.\end{equation}
\textbf{Corollary 3.3.} \emph{Let one of conditions I-III of
Theorem 3.1 be fulfilled and}\begin{equation}
\frac{t}{[a(t)]^{\alpha}}{\to}0,\quad
t{\to}\infty.\end{equation}\emph{Then the following equalities are
true:}\begin{equation} 1)\quad \lim{{p_{\alpha}(t,a(t))}}=1,\quad
t{\to}\infty.\end{equation}
\begin{equation} 2)\quad \lim {P[\sup{|X(\tau)|}}>a(t)]=0\quad
0{\leq}\tau{\leq}t,\quad t{\to}\infty.
\end{equation}
Corollary 3.3 follows from (3.3) and the relation \begin{equation}
\lim{p_{\alpha}(t)}=1 ,\quad t{\to}0 .\end{equation}
\textbf{Corollary 3.4.} \emph{Let one of conditions I-III of
Theorem 3.1 be fulfilled and}\begin{equation}
\frac{t}{[a(t)]^{\alpha}}{\to}T,\quad 0<T<\infty,\quad
t{\to}\infty.\end{equation}\emph{Then the following equality is
true:}\begin{equation}
\lim{{p_{\alpha}(t,a(t))}}=p_{\alpha}(T),\quad
t{\to}\infty.\end{equation} Corollary 3.4 follows from (3.3).
\section{Wiener Process}
We consider separately the important special case when $\alpha =2$
(Wiener process).  In this case the kernel $\Phi_{2}(x,t,a)$ of
the operator $B_{\alpha}$ coincides with Green's function  (see
[3]) of the equation \begin{equation} -\frac
{d^{2}y}{dx^{2}}=f(x),\quad -a{\leq}x{\leq}a \end{equation} with
the boundary conditions
\begin{equation}
y(-a)=y(a)=0. \end{equation}It is easy to see that
\begin{equation}
\Phi_{2}(x,t,a)=\frac{1}{2a}\begin{cases}(t+a)(a-x)
&\text{$-a{\leq}t{\leq}x$}\\ (a-t)(a+x) &\text{$x{\leq}t{\leq}a$}
\end{cases}\end{equation} According to formula (4.3) equality (3.1)
is true when $\alpha=2$ too. Hence the formula \begin{equation}
p_{2}(t,a)=p_{2}(t/a^{2},1) \end{equation} is true as well.\\ Let
$a=1.$ The eigenvalues of problem (4.1),(4.2) have the form
\begin{equation}
\mu_{n}=(n\pi/2)^{2},\quad n=1,2,3...\end{equation} The
corresponding normalized eigenfunctions are defined by the
equality
\begin{equation} g_{n}(x)=\sin{[(n\pi/2)(x-1)]}.\end{equation} The
following formula \begin{equation}
 p_{2}(t)=\sum_{j=1}^{\infty}g_{j}(0)\int_{-1}^{1}g_{j}(x)dxe^{-t\mu_{j}}
 \end{equation} holds. Formula (4.7) can be deduced in the same way
 as in the case ($0<\alpha<2,\quad \beta=0$) (see [3], [6]). Using (4.5) and
 (4.6) we rewrite (4.7) in the form \begin{equation}
 p_{2}(t)=\sum_{m=0}^{\infty}(-1)^{m}\frac{2}{\pi(m+1/2)}e^{-t[\pi(m+1/2)]^{2}}.
 \end{equation}The series (4.8) satisfies the conditions of Leibniz
 theorem. It means that  $p_{2}(t,a)$  can be
 calculated with a given precision when the parameters $t$ and $a$ are
 fixed.  From (4.4) and
(4.8) we deduce that \begin{equation} p_{2}(t,a)=
\frac{4}{\pi}e^{-t\pi^{2}/[2a(t)]^{2}}(1+o(1)),
\end{equation} where $t/[a(t)]^{2}{\to}\infty$.\\
\textbf{Proposition 4.1.} \emph{Theorem 3.1 and Corollaries
3.1-3.4 are true in case when }$\alpha=2$ too.\\ \textbf{Remark
4.1} From the probabilistic point of view  it is easy to see that
the function $p_{2}(t) , (t>0)$ is monotonic decreasing and
\begin{equation} 0<p_{2}(t){\leq}1;\quad\lim {p_{2}(t)}=1,\quad
t{\to}0.\end{equation}
\section{Iterated logarithm law, most visited sites and  first
hitting time.} It is interesting to compare our results (Theorem
3.1,
 Corollaries 3.1-3.4 and Proposition 4.1) with the well-known results mentioned in the title of
 the section.\\
1. We begin with the famous Khinchine   theorem (see[10]) about
the iterated logarithm law.\\ \textbf{Theorem 5.1.} \emph{Let
$X(t)$ be stable process ($0<\alpha<2$). Then almost surely (a.s.)
that} \begin{equation} \lim
{\frac{\sup{|X(t)|}}{(t\log{t})^{1/\alpha}|\log{|\log{t}|}|^{(1/{\alpha)+\epsilon}}}}=
\begin{cases}0 &\text{$\epsilon>0$ a.s.}\\
\infty &\text{$\epsilon=0$ a.s.}
\end{cases}\end{equation}
 We introduce the random process \begin{equation}
U(t)=\sup{|X(\tau)|}, \quad {0\leq}\tau{\leq}t \end{equation} From
Corollaries 3.1-3.4 and Proposition 4.1 we deduce the assertion.\\
\textbf{Theorem 5.2.} \emph{Let one of  conditions I-III of
Theorem 3.1 be fulfilled or let $\alpha=2$ and
\begin{equation} b(t){\to}\infty ,\quad t{\to}\infty.\end{equation} Then }
\begin{equation} b(t)U(t)/t^{1/\alpha}{\to}\infty\quad a.s.,\quad U(t)/[b(t)t^{1/\alpha}]{\to}0 \quad a.s.\end{equation}\\
In particular we have :
\begin{equation} [(\log^{\epsilon}{t})U(t)]/t^{1/\alpha}{\to}\infty\quad a.s.,\quad U(t)/[(\log^{\epsilon}{t})t^{1/\alpha}]{\to}0 \quad a.s.\end{equation}\\
when $\epsilon>0.$ We see that our approach and the classical one
have some similar points (estimation of $|X(\tau)|$), but these
approaches are essentially different. We consider the behavior of
$|X(\tau)|$ on the interval $(0,t)$, and in the classical case
$|X(\tau)|$ is considered on the interval $(t,\infty)$.\\ 2.  We
denote by $V(t)$ the most visited site of stable process $X$ up to
time $t$ (see [1]). We formulate the following result (see [1] and
references therein).\\ Let $1<\alpha<2,\quad \beta =0,\quad \gamma
>9/(\alpha-1)$. Then we have \begin{equation}
\lim{\frac{(\log{t})^{\gamma}}{t^{1/\alpha}}|V(t)|}=\infty ,\quad
t{\to}\infty\quad a.s.\end{equation}
 To this important result we add the following estimation.\\
\textbf{Theorem 5.3.} \emph{Let one of the conditions I-III of
Theorem 3.1 be fulfilled or let $\alpha=2$ and
\begin{equation} b(t){\to}\infty ,\quad t{\to}\infty.\end{equation} Then }
\begin{equation} |V(t)(t)/[b(t)t^{1/\alpha}]{\to}0 \quad a.s.\end{equation}\\
In particular we have when $\epsilon>0$:
\begin{equation}  |V(t)|/[(\log^{\epsilon}{t})t^{1/\alpha}]{\to}0 \quad a.s.\end{equation}\\
 The formulated theorem follows directly from the
inequality $U(t){\geq}|V(t)|$. \\ 3. The first hitting time
$T_{a}$ is defined by the formula
\begin{equation}T_{a}=\inf{(t{\geq}0,X(t){\geq}a)}.
\end{equation} It is obvious that
\begin{equation}P(T_{a}>t)=P[\sup{X(\tau)}<a],\quad
0{\leq}\tau{\leq}t.\end{equation} We have
\begin{equation}P(T_{a}>t){\geq}P[\sup{|X(\tau)|}<a]=p_{\alpha}(t,a),\quad
0{\leq}\tau{\leq}t.\end{equation}So our formulas for
$p_{\alpha}(t,a)$ estimate $P(T_{a}>t)$ from below.
\\\textbf{Remark 5.1.} Our results can be interpreted in terms of
the  first hitting time $T_{\pm{a}}$ one of the barriers
$x=\pm{a}$. Namely, we have
\begin{equation} P(T_{\pm{a}}>t)=p_{\alpha}(t,a).\end{equation}
The distribution of the first hitting time for the stable
processes
 is an open problem.\\
 \textbf{Remark
5.2.} Rogozin B.A. in his interesting work [5] established the law
of the overshoot distribution for the stable processes when the
existing interval is fixed.\\

\begin{center}\textbf{References}\end{center}
1. \textbf{Bass R.F., Eisenbaum N. and Shi Z.}, The Most Visited
Sites of Symmetric Stable Processes,  Probability Theory and
Related Fields, 116, (2000), 391-404.\\ 2. \textbf{Feller W.}, An
Introduction to Probability Theory and its Applications, J.Wiley
and Sons, 1971.\\ 3. \textbf{Kac M.}, On some Connections Between
Probability Theory and Differential and Integral Equations,
Proc.Sec.Berkeley Symp.Math.Stat. and Prob., Berkeley, 189-215,
1951.\\ 4. \textbf{Pozin S.M., Sakhnovich L.} Two-sided Estimation
of the Smallest Eigenvalue of an Operator Characterizing Stable
Processes, Theory Prob. Appl., 36, No.2,385-388, 1991 \\ 5.
\textbf{Rogozin B.A.,}The distribution of the first hit for stable
and asymptotically stable walks on an interval, Theory Probab.
Appl. 17, 332-338, 1972.\\ 6. \textbf{Sakhnovich L.}, Abel
Integral Equations in the theory of Stable Processes, Ukr.Math.
Journ., 36:2, 193-197 , 1984.\\ 7. \textbf{Sakhnovich L.} Integral
Equations  in the theory of Stable Processes, St.Peterburg
Math.J.,4, No.4 1993,819-829.\\ 8. \textbf{Sakhnovich L.}The
Principle of Imperceptibility of the Boundary in the Theory of
Stable Processes , St.Peterburg Math.J., 6, No.6, 1995,
1219-1228.\\ 9. \textbf{Sakhnovich L.}, Integral Equations with
Difference Kernels, Operator Theory, v.84, 1996, Birkhauser.\\ 10.
\textbf{Thomas M., Barndorff O.}(ed.) Levy Processes ; Theory and
Applications, Birkhauser, 2001.\\ 11. \textbf{Widom H.}, Stable
Processes and Integral Equations, Trans. Amer. Math. Soc. 98,
430-449, 1961.\\

\end{document}